\documentclass[11pt]{article}

\usepackage{graphicx,color}

\usepackage{amssymb}

\usepackage[english,francais]{babel}

\newtheorem{e-proposition}[theorem]{Proposition}

\newtheorem{e-definition}[theorem]{Definition\rm}

\newtheorem{theoreme}{Th\'eor\`eme}[section]
\newtheorem{lemme}[theoreme]{Lemme}
\newtheorem{proposition}[theoreme]{Proposition}

\newtheorem{remarque}{\it Remarque}

\def\og{\leavevmode\raise.3ex\hbox{$\scriptscriptstyle\langle\!\langle$~}}
\def\fg{\leavevmode\raise.3ex\hbox{~$\!\scriptscriptstyle\,\rangle\!\rangle$}}
\newcommand{\boK}{\mathcal{K}}
\newcommand{\boL}{\mathcal{L}}
\newcommand{\Z}{\mathbb{Z}}
\newcommand{\R}{\mathbb{R}}
\newcommand{\sign}{\textrm{Sign}}
\newcommand{\ba}[1]{\overline{#1}}

\begin{document}

\title{Comportement \`a l'infini du graphe gordien des n\oe{}uds
\footnote{\emph{keywords}: n\oe{}ud gordien, quasi-isom\'etrie, signature. \quad \quad \quad \emph{2000 Mathematics Subject Classification:} 57M27.}}

\author{Julien March\'e
\footnote{Institut de Math\'ematiques de Jussieu, \'Equipe ``Topologie et G\'eom\'etries Alg\'ebriques''
Case 7012, Universit\'e Paris VII, 75251 Paris CEDEX 05, France. 
\quad e-mail: \texttt{marche@math.jussieu.fr}}}
\maketitle

\selectlanguage{english}
\begin{abstract}
We study the gordian graph of all knots in $\R^3$: two knots are adjacent if they differ by a single crossing change. We prove that this graph contains isometrically an infinite countable tree with infinite valency, and that the complement of any finite subset is connected. 

\vskip 0.5\baselineskip

\selectlanguage{francais}
\begin{center}{\bf R\'esum\'e}\end{center}

On \'etudie le graphe gordien des n\oe{}uds dans $\R^3$: deux n\oe{}uds sont adjacents si on passe de l'un \`a l'autre en changeant un croisement. On prouve que ce graphe contient isom\'etriquement un arbre infini d\'enombrable de valence infinie et que le compl\'ementaire de tout sous-ensemble fini est connexe.

\end{abstract}

Soit $\boK$ l'ensemble des classes d'isotopie de n\oe{}uds orient\'es dans $S^3$.
On consid\`ere le {\it graphe gordien} dont les sommets sont les \'el\'ements de $\boK$ et tel que deux n\oe{}uds sont reli\'es par une ar\^ete si on peut passer de l'un \`a l'autre en changeant un croisement. On notera $U$ le n\oe{}ud trivial et $d_K$ la distance induite sur $\boK$ par le graphe. On appelle gordien tout n\oe{}ud $L$ v\'erifiant $d_K(L,U)=1$.

Depuis que le graphe gordien a \'et\'e introduit (voir \cite{japs}), peu de r\'esultats ont \'et\'e obtenus sur sa structure. On s'est essentiellement int\'eress\'e \`a la structure locale du graphe en montrant par exemple que chaque ar\`ete de $\boK$ fait partie d'un sous-graphe complet infini.

A contrario, J.M.~Gambaudo et \'E.~Ghys ont propos\'e de s'int\'eresser au comportement \`a l'infini du graphe, c'est-\`a-dire, \`a quasi-iom\'etrie pr\`es. Dans cette optique, ils ont prouv\'e qu'on pouvait plonger le graphe $\Z^n$ quasi-isom\'etriquement dans $\boK$ pour toute valeur de $n$ (\cite{gamb,ghys}).

R\'epondant \`a une question d'\'E. Ghys, on d\'emontre dans cet article que l'on peut plonger isom\'etriquement un arbre infini d\'enombrable de valence infinie dans le graphe gordien. Avec les m\^emes techniques, on prouve que le compl\'ementaire de tout ensemble fini dans $\boK$ est connexe.

On rappelle qu'il existe un invariant de $K$ appel\'e {\it signature} $\sigma_K:S^1\to \Z$. Cet invariant v\'erifie pour tout $z\in S^1$ et $K_1,K_2\in\boK$ l'in\'egalit\'e suivante:
$$|\sigma_{K_1}(z)-\sigma_{K_2}(z)|\le 2 d_{\boK}(K_1,K_2).$$

\section{Plongement d'un arbre infini dans le graphe gordien}

Soit $T$ l'arbre trivalent infini d\'enombrable sur lequel on a ajout\'e un sommet sur une ar\^ete appel\'e {\it racine}. De cette mani\`ere, chaque sommet a deux {\it enfants} et tous les sommets ont un unique {\it parent} sauf la racine.

\begin{theoreme}
Il existe un plongement quasi-isom\'etrique de $T$ dans $\boK$.
\end{theoreme}
\begin{remarque}
La meme preuve permet de montrer que l'on peut plonger un arbre r\'egulier de valence infinie. Ce plongement est en fait isom\'etrique si on multiplie par 2 la distance de $T$.
\end{remarque}

Le reste de cette partie consiste \`a prouver cette proposition.
La m\'ethode consiste \`a construire ce plongement r\'ecursivement \`a partir du n\oe{}ud trivial qui est associ\'e \`a la racine. Chaque autre sommet sera obtenu en faisant la somme connexe de son parent avec deux n\oe{}uds gordiens de sorte que sa signature soit ``assez compliqu\'ee''.
Pour pr\'eciser ceci, on utilise le lemme suivant qui sera prouv\'e dans la troisi\`eme section.

\begin{lemme}\label{lemmefond}
Les propri\'et\'es suivantes sont v\'erifi\'ees:
\begin{itemize}
\item
Soit $\Delta$ un polyn\^ome normalis\'e, c'est-\`a-dire un polyn\^ome $\Delta\in\Z[t,t^{-1}]$ v\'erifiant $\Delta(t^{-1})=\Delta(t)$ et $\Delta(1)=1$. Il existe alors un n\oe{}ud gordien qui admette $\Delta$ pour polyn\^ome d'Alexander.
\item
La fonction signature d'un n\oe{}ud gordien de polyn\^ome d'Alexander normalis\'e $\Delta$ est la fonction qui \`a tout nombre complexe $z$ de module 1 associe $1-\sign(\Delta(z))$ ou la fonction oppos\'ee.
\end{itemize}
\end{lemme}

Consid\'erons pour tout $p$ impair le polyn\^ome normalis\'e du n\oe{}ud torique de param\`etres $(p,2)$ d\'efini par $D_p(t)=t^{-(p-1)/2}(t^p+1)/(t+1)$.
Gr\^ace au lemme on peut choisir un n\oe{}ud $K_p$ qui est gordien et a $D_p$ pour polyn\^ome d'Alexander ($K_p$ n'est pas le n\oe{}ud torique de param\`etres $(p,2)$ car ce dernier n'est pas gordien). On suppose de plus que la signature de $K_p$ est positive (sinon, on prend son image dans un miroir).
Le dessin ci-dessous repr\'esente les endroits du cercle o\`u la signature de $K_p$ vaut 2 (elle vaut 0 sur le compl\'ementaire et 1 aux bornes). Notons $A_p$ ce sous-ensemble du cercle.
\begin{figure}[htbp]
\begin{center}
\input{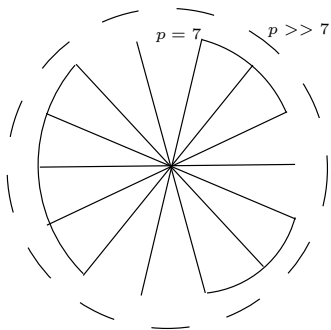}
\label{sign}
\caption{Signature des n\oe{}uds $K_p$}
\end{center}
\end{figure}

\begin{lemme}
On peut trouver une suite strictement croissante de nombres impairs $p_n$ telle que pour toute suite finie $i_1<i_2<\ldots<i_k$ et $\epsilon_1,\ldots \epsilon_k \in \{\pm 1\}$, on ait:
$$A_{p_{i_1}}^{\epsilon_1}\cap \cdots \cap A_{p_{i_k}}^{\epsilon_k}\ne \emptyset$$ o\`u $A_p^1=A_p$ et $A_p^{-1}$ est le compl\'ementaire de $A_p$.
\end{lemme}

{\bf D\'emonstration:} On construit cette suite r\'ecursivement. Prenons par exemple $p_1=3$. En choisissant $p_2$ assez petit, on peut s'assurer que chaque composante connexe de $A_{p_1}$ et $A_{p_1}^{-1}$ intersecte \`a la fois $A_{p_2}$ et $A_{p_2}^{-1}$. On continue cette construction \`a partir de $p_2$, ce qui prouve le lemme. On peut par exemple d\'efinir $p_n$ par la formule $p_n=(2n+1)(2n-1)\cdots 3$ pour $n\ge 1$. $\square$

\begin{figure}[htbp]
\begin{center}
\input{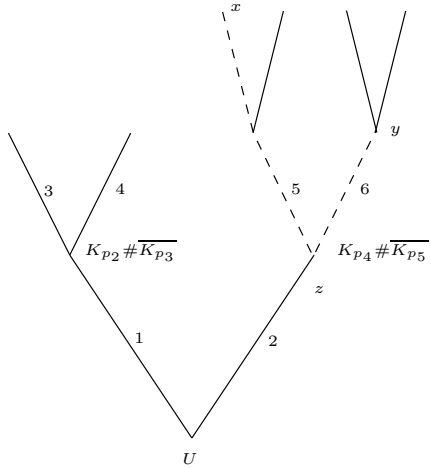}
\label{arb}
\caption{Le plongement $\phi$}
\end{center}
\end{figure}

On construit maintenant une application $\phi:T\to \boK$. On num\'erote toutes les ar\^etes de $T$. Chaque sommet $x$ de $T$ diff\'erent de la racine est reli\'e \`a son parent $y$ par une ar\^ete num\'erot\'ee $n$. Notons  $\#$ la somme connexe et $\ba{K}$ l'image miroir d'un n\oe{}ud $K$. On pose
 $$\phi(x)=\phi(y)\#K_{p_{2n}}\#\ba{K_{p_{2n+1}}}.$$
Cette construction est r\'esum\'ee dans la figure \ref{arb}.

Comme pour tout $p$ le n\oe{}ud $K_p$ est gordien, on a $d_{\boK}(\phi(x),\phi(y))\le 2$. Ainsi, on a $\forall x,y\in T$, $d_{\boK}(\phi(x),\phi(y))\le 2 d_T(x,y)$. On passe maintenant \`a l'autre in\'egalit\'e. Prenons deux \'el\'ements quelconques $x$ et $y$ de $T$ et notons $z$ leur premier parent commun. Les n\oe{}uds $\phi(x)$ et $\phi(y)$ sont obtenus \`a partir de $\phi(z)$ par une suite de sommes connexes. Pour expliciter celle-ci, supposons que

$\phi(x)=\phi(z) \#K_{p_{2n_1}}\# \ba{K_{p_{2n_1+1}}}\#\cdots \#K_{p_{2n_k}}\# \ba{K_{p_{2n_k+1}}}$ et 

$\phi(y)=\phi(z)\#K_{p_{2m_1}}\# \ba{K_{p_{2m_1+1}}}
\#\cdots \#K_{p_{2m_l}}\# \ba{K_{p_{2m_l+1}}}$ de sorte que $d_T(y,z)=k+l$.

D'apr\`es le lemme, $A_{p_{2n_1}}\cap A_{p_{2n_1+1}}^{-1}\cap\cdots A_{p_{2n_k}}\cap A_{p_{2n_k+1}}^{-1}\cap 
A_{p_{2m_1}}^{-1}\cap A_{p_{2m_1+1}}\cap\cdots \cap A_{p_{2m_l}}^{-1}\cap A_{p_{2m_l+1}}\ne \emptyset$. On choisit $\theta$ dans cet ensemble de sorte que $
\sigma_{\phi(x)}(\theta)=\sigma_{\phi(z)}(\theta)+4k$ et $
\sigma_{\phi(y)}(\theta)=\sigma_{\phi(z)}(\theta)-4l$.

L'identit\'e $|\sigma_{\phi(x)}(\theta)-\sigma_{\phi(y)}(\theta)|=4(k+l)=4 d_T(x,y)$  prouve que $d_{\boK}(\phi(x),\phi(y))\ge 2d_T(x,y)$.
L'application $\phi$ v\'erifie ainsi $\forall x,y\in T$, $d_{\boK}(\phi(x),\phi(y))= 2 d_T(x,y)$. C'est donc bien une quasi-isom\'etrie.

\section{Connexit\'e \`a l'infini}

Toujours \`a l'aide du lemme, on peut prouver le r\'esultat suivant:

\begin{proposition}
Le compl\'ementaire de tout ensemble fini de $\boK$ est connexe.
\end{proposition}

{\bf D\'emonstration:}
Soit $\{L_1,\ldots,L_n\}$ un sous-ensemble fini de $\boK$. On note alors $\Delta_i$ le polyn\^ome d'Alexander de $K_i$ et $l(\Delta_i)$ la plus petite distance entre deux racines de $\Delta_i$ cons\'ecutives sur le cercle unit\'e. Enfin, on pose $l=\min_i l(\Delta_i)$. Pour $p$ assez grand, on peut trouver comme dans la section pr\'ec\'edente un n\oe{}ud gordien $K_p$ de polyn\^ome d'Alexander $D_p$ tel que $l(D_p)<l$.

Par des simples consid\'erations de signature, on remarque que pour tout $i$ entre 1 et $n$, on a $L_i\# K_p \notin \{L_1,\ldots,L_n\}$.

On conclut alors de la fa\c{c}on suivante: soit  $\boL=\{L_1,\ldots,L_n\}$ un sous-ensemble fini de $\boK$, et $L_{\alpha}, L_{\beta}\notin \boL$. Par connexit\'e de $\boK$, il existe un chemin $L_{\alpha}=L_{\alpha_0},L_{\alpha_1},\ldots,L_{\alpha_k}=L_{\beta}$ de n\oe{}uds contigus. 

Appliquant la remarque ci-dessus \`a l'ensemble $\boL\cup\{L_{\alpha_i}\}$, on trouve un n\oe{}ud gordien $K_p$ qui v\'erifie $L_{\alpha_i}\#K_p \notin \boL$ pour tout $i$. Le chemin $L_{\alpha_0},L_{\alpha_0}\#K_p,\ldots,L_{\alpha_k}\# K_p,L_{\alpha_k}$ est un chemin qui relie $L_{\alpha}$ et $L_{\beta}$ dans le compl\'ementaire de $\boL$. $\square$

\section{D\'emonstration du lemme \ref{lemmefond}}

Ce r\'esultat utilise des techniques tr\`es classiques et, bien que non publi\'e \`a notre connaissance, il est certainement connu des sp\'ecialistes. Par souci de compl\'etude, on donne ici une preuve.
Consid\'erons le n\oe{}ud de la figure \ref{gordien}:
\begin{figure}[htbp]
\begin{center}
\input{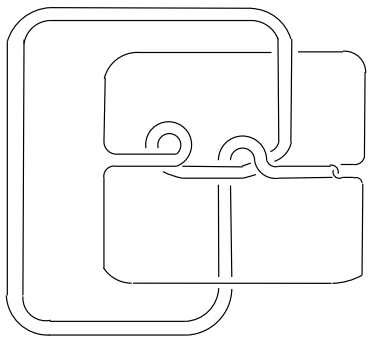}
\label{gordien}
\caption{Exemple de n\oe{}ud gordien}
\end{center}
\end{figure}

C'est clairement un n\oe{}ud gordien, et un calcul montre que son polyn\^ome d'Alexander est $-t^{-2}+3t^{-1}-3+3t-t^2$. Cette construction sugg\`ere la g\'en\'eralisation donn\'ee par la figure 4. On utilise pour les calculs la pr\'esentation chirurgicale des n\oe{}uds. C'est une mani\`ere alternative aux matrices de Seifert d'acc\'eder au polyn\^ome d'Alexander et \`a la signature des n\oe{}uds. 

Un n\oe{}ud gordien $K$ se d\'efait par changement d'un croisement. Cela signifie qu'il existe une courbe parall\'elis\'ee $L$ homologue \`a 0 dans le compl\'ementaire du n\oe{}ud trivial qui donne par chirurgie la vari\'et\'e $S^3\setminus K$.

L'interpr\'etation de la signature et du polyn\^ome d'Alexander de $K$ en fonction de la matrice d'enlacement \'equivariant des composantes de chirurgie (voir \cite{surg}) nous donne l'information cl\'e suivante: soit $P$ l'auto-enlacement \'equivariant de $L$. C'est un polyn\^ome de Laurent en $t$ qui v\'erifie $\Delta(K)=P(1).P$ et $\sigma_K(z)=\sign(P(z))-\sign(P(1))=\pm(1-\sign(\Delta(z)))$.

Il reste donc \`a \'etablir le premier point. Cela revient \`a prouver que pour tout polyn\^ome normalis\'e $P$, on peut construire une composante de chirurgie non nou\'ee et homologue \`a 0 dans le compl\'ementaire du n\oe{}ud trivial qui a $P$ pour auto-enlacement \'equivariant.

Fixons une famille d'entiers $a_0,\ldots,a_n$ et consid\'erons la courbe de la figure 4 dans laquelle on a identifi\'e la partie sup\'erieure et la partie inf\'erieure. De cette mani\`ere, la courbe est naturellement plong\'ee dans un tore plein.

\begin{figure}[htbp]
\begin{center}
\input{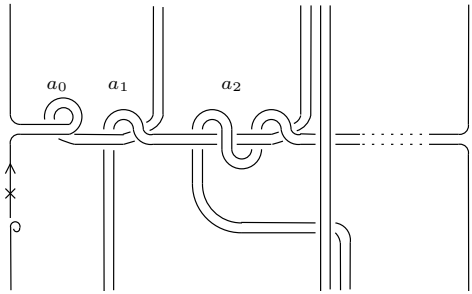}
\caption{Courbe de chirurgie}
\end{center}
\end{figure}

On calcule l'auto-enlacement \'equivariant de cette courbe en comptant le nombre de points d'intersections avec l'information relative au plongement.
Utilisant le formalisme des diagrammes de cordes, on code les croisements \`a l'aide du diagramme de la figure 5.

\begin{figure}[htbp]
\begin{center}
\input{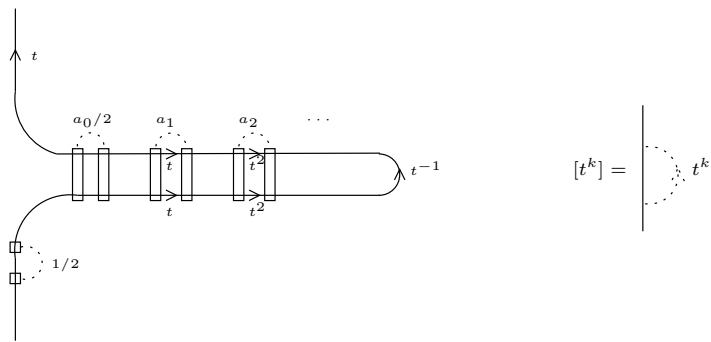}
\caption{Calcul de l'auto-enlacement \'equivariant}
\end{center}
\end{figure}
Le symbole $t$ d\'esigne le g\'en\'erateur du groupe fondamental du tore plein dirig\'e vers le haut. Les bo\^ites reli\'ees par une ligne pointill\'ee d\'esignent l'ensemble des ar\^etes devant \^etre reli\'e par un diagramme de cordes. Chaque corde relie deux points sur la courbe $L$: le coloriage de la corde d\'esigne la classe d'homotopie du cycle form\'e par la composante coup\'ee et la corde.

En d\'eveloppant ces diagrammes, on obtient le calcul d\'etaill\'e dans la figure 6.

\begin{figure}[htbp]
\begin{center}
\input{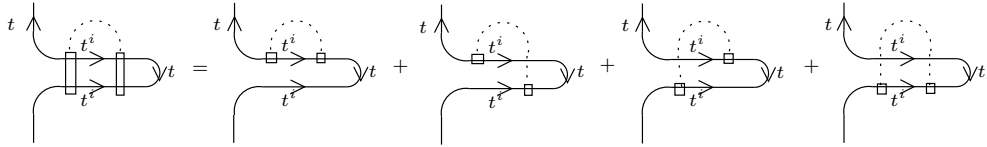}
\caption{Contributions \`a l'auto-intersection}
\end{center}
\end{figure}
Dans les notations de la figure 5, le diagramme de corde de la figure 6 est \'egal \`a l'expression $[t^i]-[t^{i}t]-[t^{-i}t]+[t^{-i}]=[(t^i+t^{-i})(1-t)]$.

Au final, le diagramme de corde associ\'e \`a la figure 5 est
$\frac{1}{2}[1]+\frac{a_0}{2}[2(1-t)]+\sum_{i\ge 1}a_i[(t^i+t^{-i})(1-t)]$.

On obtient le polyn\^ome d'Alexander en prenant la partie sym\'etrique de cette expression \`a savoir $1+a_0(2-t-t^{-1})+\sum_{i\ge 1} a_i (t^i+t^{-i})(2-t-t^{-1})$. Or, il est clair que tout polyn\^ome normalis\'e s'\'ecrit de mani\`ere unique sous cette forme pour un choix convenable des coefficients $a_i$. Cela prouve bien que tout polyn\^ome normalis\'e est le polyn\^ome d'Alexander d'un n\oe{}ud gordien.


\section*{Remerciements}
Je tiens \`a remercier \'E.~Ghys pour avoir initi\'e et suivi cette \'etude, ainsi que G.~Masbaum et P.~Vogel pour leurs remarques utiles.

\end{document}